\documentclass[11pt]{article}  
\usepackage{times,amsmath,amsfonts,amsthm,amssymb,eucal}
\usepackage{array}
\usepackage{color}
\usepackage{enumerate}
\usepackage{epsfig}
\usepackage{mathrsfs}
\usepackage{fancyhdr}
\usepackage{graphicx}
\usepackage{caption}
\usepackage{subcaption}
\usepackage{graphicx}
\usepackage{makeidx}
\usepackage{multirow}
\newcommand{\E}{\mathrm{E}}
\newcommand{\Expect}{{\rm I\kern-.3em E}}
\newtheorem{theorem}{Theorem}

\newtheorem{proposition}{Proposition}
\newtheorem{lemma}{Lemma}
\newtheorem{corollary}{Corollary}
\newtheorem{definition}{Definition}
\newtheorem{assumption}{Assumption}
\def\Theorem{\begin{theorem}\sl}
\def\EndTheorem{\end{theorem}}
\def\Proposition{\begin{proposition}\sl}
\def\EndProposition{\end{proposition}}
\def\Lemma{\begin{lemma}\sl}
\def\EndLemma{\end{lemma}}
\def\Corollary{\begin{corollary}\sl}
\def\EndCorollary{\end{corollary}}
\def\Definition{\begin{definition}\sl}
\def\EndDefinition{\end{definition}}
\numberwithin{equation}{section}
\begin{document}
\title{ \textbf{A Note on Bootstrapping M-estimates from Unstable AR($2$) Process with Infinite Variance Innovations}}   
\author{Maryam Sohrabi and Mahmoud Zarepour\\ University of Ottawa, Ontario, Canada}
\date{\today}
\maketitle
\begin{abstract}
The limiting distribution for M-estimates in a non-stationary autoregressive model with heavy-tailed error is computationally intractable. To make inferences based on the  M-estimates,  the bootstrap procedure can be used to approximate the sampling distribution. In this paper, we show that the bootstrap scheme with $m=o(n)$ resampling sample size when $m/n \to 0$  is  approximately valid in a multiple unit roots  time series with innovations in the domain of attraction of a stable law with index $0<\alpha\leq2$.\\

\vspace{9pt}
\end{abstract}
\section {Introduction}
\label{intro}
Consider the autoregressive process of order $p$ (AR($p$))
\begin{align}
\label{1}
 \phi(B)X_t=\epsilon_{t},
\end{align}
where $B$ is the backward operator and
\begin{align}
\label{12}
 \phi(z)=1-\phi_1z-\phi_2z^2-\cdots-\phi_pz^p.
\end{align}
The errors $\{\epsilon_{t}\}$ in \eqref{1} form a sequence of independent and identically distributed (i.i.d.) random variables in the domain of attraction of a symmetric stable law with index $0<\alpha\le2$. For  $0<\alpha<2$, this is equivalent to the following assumptions,
$$ P(|\epsilon_{1}|>x)=x^{-\alpha}L(x)$$
for some  slowly varying  function $L$ at $\infty$ with $\alpha>0$ and
\begin{align*}
\frac{P(\epsilon_{1}>x)}{P(|\epsilon_{1}|>x)}\to p
\quad \text{and} \quad
\frac{P(\epsilon_{1}\le -x)}{P(|\epsilon_{1}|>x)}\to q,
\end{align*}
as $x\to\infty$, $0\le p\le 1$ and $q=1-p$.  If $\{\epsilon_{t}\}$ has a symmetric distribution then $p=q=1/2$.

The model \eqref{1} is referred to as a non-stationary  autoregressive  time series, if the characteristic polynomial $\phi(\cdot)$ has at least one root on  the boundary of the unit circle. The problem to conduct asymptotic inference for time series with unit roots has been a challenging topic of interest.  The extensive studies about autoregressive time series models under the heavy-tailed hypothesis  include  Knight (1989), Chan and Tran (1989), Phillips (1990),  Davis, Knigh, and Liu (1992), Davis and Wu (1997), Tanaka (2008), Samarakoon and Knight (2009), Moreno and Romo (2012), and Chan and Zhang (2012).

When the  residuals have regularly varying tail probabilities, least square (LS) estimation methods can exhibit rather poor power performance. Thus,  it is important to consider estimation and inference procedures which are robust to departures from finite variance condition. One way to achieve robustness is the use of M-estimate method.
For a given loss function $\rho(x)$, the M-estimate $\hat{\Phi}=(\hat{\phi}_{1},\ldots,\hat{\phi}_{p})$ of $\Phi=({\phi}_{1},\ldots,{\phi}_{p})$ minimizes
the objective function

\begin{align*}
\sum_{t=p+1}^{n}\rho(X_t-\beta_{1}X_{t-1}-\ldots-\beta_{p}X_{t-p}),
\end{align*}
with respect to $\left(\beta_1 , \ldots,\beta_p\right)$, where $\rho$ is an almost everywhere differentiable convex function. This guarantees the uniqueness of the solution. For more details see Davis et al. (1992). Knight (1989) finds the asymptotic distribution of M-estimators of the autoregressive parameter of an infinite-variance random walk. The results establish that self-normalized M-estimates are asymptotically normal and their rate of convergence is higher than the LS estimates. Davis et al. (1992)  study the problem of estimating autoregressive parameters when the observations are from a stationary  AR($p$)
process with innovations in the domain of attraction of a stable law. Sohrabi (2016)  present the asymptotic distribution of M-estimators for parameters in unstable AR($p$) processes when the innovations are assumed to be in the domain of attraction of a symmetric stable law with index $0<\alpha\le2$.

 While the asymptotic theory for M-estimates is well understood, the limiting distributions are generally intractable. This prevents the use of the asymptotic distribution for inference purposes such as for the construction of confidence intervals. Davis and Wu (1997) consider bootstrapping M-estimates for the  stationary AR($p$) processes. Moreno and Romo (2012) use bootstrap resamples to estimate the percentiles of the limiting distribution of M-estimates in a random walk process with errors in the domain of attraction of a law.   In this paper, we
investigate the bootstrap for approximating the distribution of M-estimates in the unstable autoregressive processes. Due the complexity of  the limiting distribution for AR($p$) processes given in Sohrabi (2016), we only show that the subsampling bootstrap procedure is asymptotically valid for  M-estimates in an unstable AR(2) process.
The rest of the article is organized as follows. In Section 2, we set up the problem and the limiting distribution of M-estimates is presented for an AR(2) process with repeated unit roots. The bootstrapping M-estimates for the unstable AR(2) process when the errors belong to the domain of attraction of a stable law with index $0<\alpha\leq2$ are considered in Section 3.

\section{Asymptotic Theory for M-Estimates in Unstable AR(2) Processes with Infinite Variance Innovations}
In this Section, for clarity and simplicity, we present the asymptotic distribution of M-estimators for parameters in a non-stationary AR(2) process when it has two real unit roots. The results for other unstable cases can be found in Sohrabi (2016).  Consider the model
\begin{align}
\label{R1}
(1-B)^2X_{t}=\epsilon_{t},
\end{align}
which is equivalent to
\begin{align*}
X_{t}=\sum_{k=1}^{t}\sum_{i=1}^{k}\epsilon_{i}.
\end{align*}
Therefore, the following property holds when we have root 1 with the multiplicity of 2:
\begin{align*}
X_{t}-X_{t-1}=\sum_{k=1}^{t}\epsilon_{k}.
\end{align*}
We assume that the innovations $\{\epsilon_t\}$ satisfy:
  \begin{assumption} (A1)
  The innovations $\{\epsilon_{t}\}$ are i.i.d. random variables in the domain of attraction of a stable law with index $0<\alpha\leq2$.
\end{assumption}
  It is common to assume symmetry for the innovations. Note that for $0< \alpha< 1$, symmetry is not required. Therefore, for $1\leq\alpha\leq2$ we impose symmetry on the innovations; i.e., $p=q=1/2$. However, for $1<\alpha\leq2$ symmetry is not required, if $\E(\epsilon_1)=0$.

\noindent Assumption A1 implies that:
\begin{align}
\label{4}
S_{n}(t)= a_{n}^{-1}\sum_{k=1}^{[nt]}\epsilon_k\ \overset{d}\to \ S(t) \quad \mbox{in} \quad D[0,1],
\end{align}
where $\{a_n\}$  is a sequence of positive constants defined as  $a_n=\inf\{x:P[|X_1|>x]\leq n^{-1}\}$ and $\overset{d}\to $ denotes here convergence in distribution with respect to the Skorohod topology. Moreover, $S(\cdot)$ is a stable process and its representation is as follows:
\begin{align}
\label{6}
S\left( t \right) =
\left\{
{\begin{array}{*{20}{l}}
{\sum\nolimits_{k = 1}^\infty  {{\delta_k}\Gamma _k^{-1/\alpha }I\left( {{U_k} \leq t} \right)} }&{\rm{if}\,0 < \alpha  < 2,}\\
{\text{standard Brownian motion}}&{\rm{if}\,\alpha  = 2,}
\end{array}} \right.
\end{align}
where $\{U_{k}\}$ is a sequence of i.i.d.  $U[0,1]$ random variables and $\{\delta_k\}$ is a sequence of i.i.d. random variables such that $P(\delta_k=1)=p$, $P(\delta_k=-1)=q$, and $p+q=1$. Also, $\Gamma_1,\Gamma_2,\ldots$ are the arrival times of a Poisson process with Lebesgue mean measure and $\{U_{k},\Gamma_{k},\delta_{k}\}$ are mutually independent. Moreover, we impose the  following assumptions on the function $\rho(\cdot)$.
\begin{assumption}
(A2)  Let $\rho$ be a convex and twice differentiable function, and take $\psi=\rho'$.
\end{assumption}
\begin{assumption}
(A3) $\E(\psi(\epsilon_1))=0$ and $\E(\psi^{2}(\epsilon_1))<\infty$.
\end{assumption}
\begin{assumption} (A4)
$0<|\E(\psi'(\epsilon_{1}))|<\infty$  and $\psi'(\cdot)$ satisfies the Lipschitz- continuity condition; i.e., there exists a real constant $k\geqslant 0$ such that for all $x$ and $y$,
 \begin{align*}
|\psi'(x)-\psi'(y)|\leq k|x-y|.
\end{align*}
\end{assumption}
\noindent Note that for the Assumptions A2-A4, sometimes $\rho'$ does not exist everywhere. In this case, although $\rho'$ is not differentiable at a countable number of points, the results will usually hold with some additional complexity in the proofs. To derive the main result of this paper, we assume that conditions A1-A4 hold and we  define the following process on the Skorohod space $D[0,1]$:
{\small
\begin{align}
\label{3}
W_{n}(t)   & = n^{-1/2}\sum_{k=1}^{[nt]}\psi(\epsilon_k),
\end{align}}
where $[x]$ stands for integer part of $x$. It is well known that $W_n(\cdot)\overset{d}\to W(\cdot)$, a standard Brownian-motion process. Similar to Theorem 4 of Resnick and Greenwood (1979), we can show that
\begin{align*}
\left(
{\begin{array}{*{20}{l}}
S_{n}(\cdot)\\
W_{n}(\cdot)
\end{array}} \right)
 \overset{d}\to
 \left(
{\begin{array}{*{20}{l}}
S(\cdot)\\
W(\cdot)
\end{array}} \right)
\end{align*}
on $D[0,1]\times D[0,1]$, where $S(\cdot)$ and $W(\cdot)$ are independent. We define the following process
\begin{align*}
A_n(u,v) = \sum_{t=3}^n\left[\rho\left({\epsilon_t}-n^{-1/2}a_{n}^{-1}u\left(X_{t-1}-X_{t-2}\right)
-n^{-3/2}a_{n}^{-1}vX_{t-2}\right)-\rho(\epsilon_t)\right],
\end{align*}
where $(u,v)=\left(n^{1/2}a_n(\hat{\phi}_{1}-2),n^{3/2}a_n\left((\hat{\phi}_{1}-2)+(\hat{\phi}_{2}+1)\right)\right)^T$ is the minimizer of $A_n(u,v)$. Using the Taylor series expansion of each summand of $A_n$ around $u=0$ and $v=0$, we get
\begin{align}
\label{i4}
\nonumber A_n(u,v)
&=-un^{-1/2}a_n^{-1} \sum_{t=3}^n\left(X_{t-1}-X_{t-2}\right)\psi(\epsilon_t)\\
\nonumber&-vn^{-3/2}a_n^{-1} \sum_{t=3}^nX_{t-2}\psi(\epsilon_t)\\
\nonumber&+\frac{1}{2}u^2n^{-1}a_n^{-2}\sum_{t=3}^n\left(X_{t-1}-X_{t-2}\right)^2\psi'(c_t^{n})\\
\nonumber&+ \frac{1}{2}v^2n^{-3}a_n^{-2}\sum_{t=3}^nX_{t-2}^{2}\psi'(c_t^{n})\\
&+ uvn^{-2}a_n^{-2}\sum_{t=3}^n{X_{t-2}}\left(X_{t-1}-X_{t-2}\right)\psi'(c_t^{n}).
\end{align}
where ${c_t}^{n}$ lies between $\epsilon_t$ and ${\epsilon_t}-n^{-1/2}a_{n}^{-1}u\left(X_{t-1}-X_{t-2}\right)-n^{-3/2}a_{n}^{-1}vX_{t-2}$. Using the fact that $\psi'$ is Lipschitz-continuous:
\begin{align*}
|\psi'(\epsilon_t)-\psi'({c_t}^{n})|\le \lambda\big{|}n^{-1/2}a_{n}^{-1}u\left(X_{t-1}
-X_{t-2}\right)+n^{-3/2}a_{n}^{-1}vX_{t-2} \big{|}.
\end{align*}
Asymptotically, $\psi'({c_t}^{n})$ can be replaced by  $\psi'(\epsilon_t)$ in \eqref{i4}. For simplicity we only consider the third term of $A_n$ (the proof is similar for the other terms in \eqref{i4}). We have
\begin{align*}
u^{2}n^{-1}a_{n}^{-2}&\sum_{t=3}^n\left(X_{t-1}-X_{t-2}\right)^2|\psi'(\epsilon_t)-\psi'({c_{t}}^{n})|\\
&\le k u^{2}n^{-1}a_{n}^{-2}\sum_{t=3}^n\left(X_{t-1}-X_{t-2}\right)^2|n^{-1/2}a_{n}^{-1}u\left(X_{t-1}
-X_{t-2}\right)+n^{-3/2}a_{n}^{-1}vX_{t-2}|\\
&\le k u^{3}n^{-1/2}m^{-1}a_{n}^{-3}\sum_{t=3}^{n} |\left(X_{t-1}-X_{t-2}\right)|^{3}\\
&+ k u^{2}vn^{-3/2}n^{-1}a_{n}^{-3}\sum_{t=3}^{n} |\left(X_{t-1}-X_{t-2}\right)|^{2}|X_{t-2}|  \overset{P}\to  0.
\end{align*}
 Furthermore, asymptotically each $\psi'(\epsilon_t)$  can be replaced by $\E\left(\psi'(\epsilon_t)\right)$ in \eqref{i4}. To appreciate why, again we consider the third term of $A_n$ where
\begin{align*}
\sum_{t=3}^n\left(X_{t-1}-X_{t-2}\right)^2\psi'(\epsilon_t)=\sum_{t=3}^n\left(X_{t-1}-X_{t-2}^*\right)^2\left[ \psi'(\epsilon_t)- {\E}\left(\psi'(\epsilon_t)\right)+ \E\left(\psi'(\epsilon_t)\right)\right].
\end{align*}
Therefore,
\begin{align}
\label{i2t}
\nonumber \frac{1}{2}n^{-1}a_n^{-2}\sum_{t=3}^n\left(X_{t-1}-X_{t-2}\right)^2\psi'({\epsilon_t}^{n})
&= n^{-1}a_n^{-2}\sum_{t=3}^n\left(X_{t-1}-X_{t-2}\right)^2\left[\psi'(\epsilon_t)- {\E}\left(\psi'(\epsilon_t)\right)\right]\\
&+  \E\left(\psi'(\epsilon_1)\right)n^{-1}a_n^{-2}\sum_{t=3}^n\left(X_{t-1}-X_{t-2}\right)^2.
\end{align}
Note that the first term of the right hand side in \eqref{i2t} approaches to zero as $n\to \infty$ since
\begin{align}
\nonumber n^{-1/2}a_n^{-2}\sum_{t=3}^n\left(X_{t-1}-X_{t-2}\right)^2\left[\psi'(\epsilon_t)- {\E}\left(\psi'(\epsilon_t)\right)\right]
\nonumber&=n^{-1/2}\sum_{t=3}^nS_{n}\left(\frac{t-1}{n}\right)\left[\psi'(\epsilon_t)- {\E}\left(\psi'(\epsilon_t)\right)\right].
\end{align}
Consequently, we have
\begin{align*}
n^{-1}a_n^{-2}\sum_{t=3}^n\left(X_{t-1}-X_{t-2}\right)^2\left[\psi'(\epsilon_t)- {\E}\left(\psi'(\epsilon_t)\right)\right]\overset{p}\to 0.
\end{align*}
 By the same justification, we can show that the preceding result holds for the last two terms of \eqref{i4}.
 Thus, the finite-dimensional distributions of $A_n$ converge weakly to those of $A$ where
\begin{align}
\label{AD}
\nonumber A(u,v) &=-u{\E}^{1/2}\left(\psi^2(\epsilon_1)\right)\int_{0}^{1}S(t)dW(t)\\
\nonumber&- v{\E}^{1/2}\left(\psi^2(\epsilon_1)\right)\int_{0}^{1}\int_{0}^{t}S(s)dsdW(t)\\
\nonumber&+ \frac{u^2}{2}{\E}\left(\psi'(\epsilon_1)\right)\int_{0}^{1}S^2(t)dt\\
\nonumber&+  \frac{v^2}{2}{\E}\left(\psi'(\epsilon_1)\right)\int_{0}^{1}\left(\int_{0}^{t}S(s)ds\right)^2dt\\
&+ uv{\E}\left(\psi'(\epsilon_1)\right)\int_{0}^{1}S(t)\int_{0}^{t}S(s)dsdt.
\end{align}
By setting the derivative of $A(u,v)$ to 0 and solving for $u$ and $v$, we have
\begin{align*}
\small
\left(
{\begin{array}{*{20}{c}}
{{n^{1/2}}{a_n}({\hat{ \phi }_{1}} -2)}\\\\
{{n^{3/2}}{a_n}({\hat{ \phi }_{1}} - 2)}+{n^{3/2}}{a_n}({\hat{ \phi }_{2}}+1)
\end{array}}
\right)
\ \overset{d}{\to}\
\Gamma _2^{ - 1}
\left(
{\begin{array}{*{20}{c}}
{\frac{{{\E^{1/2}}\left( {\psi ^2}({\epsilon _1})\right) \int_0^1 {S\left( t \right)\,dW\left( t \right)} }}{{\E\left( \psi '({\epsilon _1})\right) }}}\\\\
{\frac{{\E^{1/2}\left( {\psi ^2}({\epsilon _1})\right) \int_0^1 {\int_0^t {S\left( s \right)\,ds} \,dW\left( t \right)} }}{{\E\left( \psi '({\epsilon _1})\right) }}}
\end{array}}
\right),
\end{align*}
where
\begin{align}
\label{GD}
{\Gamma _2} =
\left(
{\begin{array}{*{20}{c}}
{\int_0^1 {{S^2}\left( t \right)\,dt}}&{\int_0^1 S\left(t\right){\int_0^t {S\left( s \right)\,ds} }\,dt }\\
{\int_0^1 {S\left( t \right)\int_0^t {S\left( s \right)\,ds} \,dt} }&{ \int_0^1 {{{\left( {\int_0^t {S\left( s \right)\,ds} } \right)}^2}\,dt}  }
\end{array}}
\right).
\end{align}

\section{Bootstrap Procedure}
The asymptotic distributions for the M-estimators obtained in the previous section are  not easily computationally tractable. This complexity grows when several real and complex unit roots appear in the time series models. Due to the complexity of the  limiting distributions, to make inferences based on M-estimates, one may consider a resampling scheme.\\

\noindent For the model defined in \eqref{R1}, consider the following steps:
\begin{enumerate}[(i)]
\item  Estimate $\phi_1$ and $\phi_2$ by $\hat{\phi}_1$ and $\hat{\phi}_2$ using M-estimate method and calculate the residuals
$$e_t=X_{t}-\hat{\phi}_1X_{t-1}-\hat{\phi}_2X_{t-2}.$$
\item  Define $\hat{F}_n(\cdot)=\frac{1}{n-2}\sum_{i=3}^{n}\varepsilon_{(e_{i}-\bar{e}\leq\cdot)}$ as the empirical distribution function of residuals , where $\bar{e}=\frac{1}{n-2}\sum_{i=3}^{n}e_{i}$ and take a sample of size $m$,  $e_1^*,\ldots,e_m^*$, from $\hat{F}_n$.
\item  The bootstrap sample $\{X_t^*\}$
 is then recursively obtained from the model
 $$X_t^*=\hat{\phi}_1X_{t-1}^*+\hat{\phi}_2X_{t-2}^*+e_t^*, \ \ t=3,\dots,m.$$
 \end{enumerate}
 Once we obtain the bootstrap series, we calculate the bootstrap estimates  from
 $$(\hat{\phi}_1^*,\hat{\phi}_2^*)=\underset{(\phi_1,\phi_2)}{\text{arg min}}\sum_{t=3}^m\rho\left(X_t^*-\phi_1X_{t-1}^*-\phi_2X_{t-2}^*\right).$$
 The following lemma is needed to derive the limiting distribution of the bootstrap estimates.
\begin{lemma}
 \label{SRS}
 Let $\{e_1^*,\ldots,e_m^*\}$ be an i.i.d. sample from $\hat{F}_n$ and $\E^*$ denotes the expectation under $\hat{F}_n$. Also, under condition A1-A4 and  with the subsampling of size $m$ such that $m/n\rightarrow 0$  as $n\rightarrow \infty$, we have
\begin{enumerate}[(i)]
\item $S_m^*(\cdot)=a_{m}^{-1}\sum_{i=1}^{[m\cdot]}e_i^*\ \overset{d}\to\  S(\cdot)$, in probability, where $S(\cdot)$ is the stable process defined in \eqref{6},
\item $\E^*(\psi(e_{1}^*))=0$,
\item $m^{-1/2}\hat{\sigma}^{-1}\sum_{i=1}^{[m\cdot]}\psi(e_{i}^*)\overset{d}\to W(\cdot), \ \text{in probability}$,
where $W(\cdot)$ is a standard Brownian motion independent of $S(\cdot)$ and $\hat{\sigma}^2=\E^*(\psi^2(e_{1}^*))=\frac{1}{n-2}\sum_{i=3}^{n}\psi^2(e_{i})\overset{p}\to \E(\psi^2(\epsilon_{1}))$,
\item $\E^*(\psi^{'}(e_{1}^*))\overset{p}\to \E(\psi^{'}(\epsilon_{1}))$.
\end{enumerate}
\end{lemma}
\noindent {\textbf{Proof.}} To prove part (i), it is enough to show that for any $k\geq 1$
\begin{align}
\label{BE}
\sum_{t=1}^m\varepsilon_{(e_t^*,a_m^{-1}\mathbf{e_{t-1}^*})}\overset{d}\to\sum_{i=1}^k\sum_{j=1}^{\infty}\varepsilon_{(\epsilon_{i,j}, {{\delta_j}\Gamma _j^{-1/\alpha }\mathbf{c_i}})}
\end{align}
in probability, where $\mathbf{e_{t-1}^*}=\left(e_{t-1}^*,\ldots,e_{t-k}^*\right)$ and $\mathbf{c_i}$ is the basis element of $\mathbb{R}^k$ with $i$th component equal to one and the rest 0. Note that $\{\delta_j\}$ and $\{\Gamma _j\}$ are as specified  in \eqref{6}. The proof of \eqref{BE}
 is  quite similar to the arguments used for Lemma 5 in Davis and Wu (1997).  Results (ii)-(iv) follow by arguments similar to those of  Proposition 4 of  Moreno and Romo (2012).  The technical details are omitted. See also  Davis and Resnick (1985) and Arcones and Gin\'{e} (1989).

\bigskip
Now, we are ready to establish the bootstrap weak convergence in probability.  Define the following process
$$ A_m^*(u,v)=\sum_{t=3}^m\left[\rho\left({e_t^*}-m^{-1/2}a_{m}^{-1}u\left(X_{t-1}^*
-X_{t-2}^*\right)-m^{-3/2}a_{m}^{-1}vX_{t-2}^*\right)-\rho(e_t^*)\right].$$
If we show that $A_m^*(u,v)\overset{d}\to A(u,v)$ in probability, where $A(u,v)$ is defined in \eqref{AD}, then the minimizer of $A_m^*(u,v)$, $\left({{m^{1/2}}{a_m}({\hat{ \phi }_{1}}^* -\hat{ \phi }_1)},{{m^{3/2}}{a_m}({\hat{ \phi }_{1}}^* - \hat{ \phi }_1)}+{m^{3/2}}{a_m}({\hat{ \phi }_{2}}^*-\hat{ \phi }_2)\right)$,  converges to the minimizer of $A(u,v)$
due to the convexity of $ A_m^*$; see  Knight (1989). To prove that $A_m^*(u,v)\overset{d}\to A(u,v)$ in probability, when $m\to\infty$, we will use a Taylor expansion of $\rho$ around $(u,v)=(0,0)$:
\begin{align}
\label{ii4}
\nonumber A_m^*(u,v)&=-um^{-1/2}a_m^{-1} \sum_{t=3}^m\left(X_{t-1}^*-X_{t-2}^*\right)\psi(e_t^*)\\
\nonumber&-vm^{-3/2}a_m^{-1} \sum_{t=3}^mX_{t-2}^*\psi(e_t^*)\\
\nonumber&+\frac{1}{2}u^2m^{-1}a_m^{-2}\sum_{t=3}^m\left(X_{t-1}^*-X_{t-2}^*\right)^2\psi'({c_t^*}^{m})\\
\nonumber&+ \frac{1}{2}v^2m^{-3}a_m^{-2}\sum_{t=3}^m{X^*_{t-2}}^{2}\psi'({c_t^*}^{m})\\
&+ uvm^{-2}a_m^{-2}\sum_{t=3}^m{X_{t-2}^*}\left(X_{t-1}^*-X_{t-2}^*\right)\psi'({c_t^*}^{m}),
\end{align}
where ${c_t^*}^{m}$ lies between $e_t^*$ and ${e_t^*}-m^{-1/2}a_{m}^{-1}u\left(X_{t-1}^*
-X_{t-2}^*\right)-m^{-3/2}a_{m}^{-1}vX_{t-2}^*$. Using the fact that $\psi'$ is Lipschitz-continuous:
\begin{align*}
|\psi'(e_t^*)-\psi'({c_t^*}^{m})|\le \lambda\big{|}m^{-1/2}a_{m}^{-1}u\left(X_{t-1}^*
-X_{t-2}^*\right)+m^{-3/2}a_{m}^{-1}vX_{t-2}^* \big{|}.
\end{align*}
Asymptotically, $\psi'({c_t^*}^{m})$ can be replaced by  $\psi'(e_t^*)$ in \eqref{ii4}. For simplicity we only consider the third term of $A_m^*$ (the proof is similar for the other terms in \eqref{ii4}). We have
\begin{align*}
u^{2}m^{-1}a_{m}^{-2}&\sum_{t=3}^m\left(X_{t-1}^*-X_{t-2}^*\right)^2|\psi'(e_t^*)-\psi'({c_{t}^*}^{m})|\\
&\le k u^{2}m^{-1}a_{m}^{-2}\sum_{t=3}^m\left(X_{t-1}^*-X_{t-2}^*\right)^2|m^{-1/2}a_{m}^{-1}u\left(X_{t-1}^*
-X_{t-2}^*\right)+m^{-3/2}a_{m}^{-1}vX_{t-2}^*|\\
&\le k u^{3}m^{-1/2}m^{-1}a_{m}^{-3}\sum_{t=3}^{m} |\left(X_{t-1}^*-X_{t-2}^*\right)|^{3}\\
&+ k u^{2}vm^{-3/2}m^{-1}a_{m}^{-3}\sum_{t=3}^{m} |\left(X_{t-1}^*-X_{t-2}^*\right)|^{2}|X_{t-2}^*|  \overset{P}\to  0,
\end{align*}
in probability, when $m\to\infty$.  Furthermore, asymptotically each $\psi'(e_t^*)$  can be replaced by $\E^*\left(\psi'(e_t^*)\right)$ in \eqref{ii4}. To appreciate why, again we consider the third term of $A_m^*$ where
\begin{align*}
\sum_{t=3}^m\left(X_{t-1}^*-X_{t-2}^*\right)^2\psi'(e_t^*)=\sum_{t=3}^m\left(X_{t-1}^*-X_{t-2}^*\right)^2\left[ \psi'(e_t^*)- {\E}^*\left(\psi'(e_t^*)\right)+ {\E}^*\left(\psi'(e_t^*)\right)\right].
\end{align*}
Therefore,
\begin{align}
\label{i2terms}
\nonumber \frac{1}{2}m^{-1}a_m^{-2}\sum_{t=3}^m\left(X_{t-1}^*-X_{t-2}^*\right)^2\psi'({e_t^*})
&= m^{-1}a_m^{-2}\sum_{t=3}^m\left(X_{t-1}^*-X_{t-2}^*\right)^2\left[\psi'(e_t^*)- {\E}^*\left(\psi'(e_t^*)\right)\right]\\
&+  {\E}^*\left(\psi'(e_1^*)\right)m^{-1}a_m^{-2}\sum_{t=3}^m\left(X_{t-1}^*-X_{t-2}^*\right)^2.
\end{align}
Note that the first term of the right hand side in \eqref{i2terms} approaches to zero as $n\to \infty$ since
\begin{align}
\label{Ai1}
\nonumber m^{-1/2}a_m^{-2}\sum_{t=3}^m\left(X_{t-1}^*-X_{t-2}^*\right)^2\left[\psi'(e_t^*)- {\E}^*\left(\psi'(e_t^*)\right)\right]
\nonumber&=m^{-1/2}\sum_{t=3}^mS_{m}^*\left(\frac{t-1}{m}\right)\left[\psi'(e_t^*)- {\E}^*\left(\psi'(e_t^*)\right)\right].
\end{align}
 Consequently, we have
\begin{align*}
m^{-1}a_m^{-2}\sum_{t=3}^m\left(X_{t-1}^*-X_{t-2}^*\right)^2\left[\psi'(e_t^*)- {\E}^*\left(\psi'(e_t^*)\right)\right]\overset{p}\to 0,
\end{align*}
in probability. By the same justification, we can show that the preceding result holds for the last two terms of \eqref{ii4}. Thus, we have
\begin{align*}
\nonumber A_m^*(u,v)&=-um^{-1/2}a_m^{-1} \sum_{t=3}^m\left(X_{t-1}^*-X_{t-2}^*\right)\psi(e_t^*)\\
\nonumber&-vm^{-3/2}a_m^{-1} \sum_{t=3}^mX_{t-2}^*\psi(e_t^*)\\
\nonumber&+\E^*\left(\psi'(e_1^*)\right)\frac{1}{2}u^2m^{-1}a_m^{-2}\sum_{t=3}^m\left(X_{t-1}^*-X_{t-2}^*\right)^2\\
\nonumber&+ \E^*\left(\psi'(e_1^*)\right)\frac{1}{2}v^2m^{-3}a_m^{-2}\sum_{t=3}^m{X^*_{t-2}}^{2}\\
&+\E^*\left(\psi'(e_1^*)\right)uvm^{-2}a_m^{-2}\sum_{t=3}^m{X_{t-2}^*}\left(X_{t-1}^*-X_{t-2}^*\right).
\end{align*}
The following two moment convergence will be necessary:
\begin{align*}
\E^*\left(\psi^2(e_1^*)\right)&=
\frac {1}{2(n-2)}\sum_{t=3}^n\left(\psi^2(e_t)+\psi^2(-e_t)\right)\\
&=\frac {1}{n-2}\sum_{t=3}^n\psi^2(e_t)=\frac {1}{n-2}\sum_{t=3}^n\psi^2\left(\epsilon_{t}+(2-\hat{\phi}_1)X_{t-1}+(-1-\hat{\phi}_2)X_{t-2}\right)\\
&= \frac {1}{n-2}\sum_{t=3}^n\psi^2(\epsilon_{t})+o_P(1)\overset{P}\to \E\left(\psi(\epsilon_1)\right),
\end{align*}
when $n\to\infty$, by a weak law of large numbers, and
\begin{align*}
\E^*\left(\psi'(e_1^*)\right)
&= \frac 1 n\sum_{t=3}^n\psi'(\epsilon_{t})+o_P(1)\overset{P}\to \E\left(\psi'(\epsilon_1)\right),
\end{align*}
when $n\to\infty$, using, again, a weak law of large numbers.
 Therefore, by applying Lemma \ref{SRS}, we  obtain the following result

\begin{align*}
\small
\left(
{\begin{array}{*{20}{c}}
{{m^{1/2}}{a_n}({\hat{ \phi }_{1}}^* -\hat{ \phi }_1)}\\\\
{{m^{3/2}}{a_m}({\hat{ \phi }_{1}}^* - \hat{ \phi }_1)}+{m^{3/2}}{a_m}({\hat{ \phi }_{2}}^*-\hat{ \phi }_2)
\end{array}}
\right)
\ \overset{d}{\to}\
\Gamma _2^{ - 1}
\left(
{\begin{array}{*{20}{c}}
{\frac{{{\E^{1/2}}\left( {\psi ^2}({\epsilon _1})\right) \int_0^1 {S\left( t \right)\,dW\left( t \right)} }}{{\E\left( \psi '({\epsilon _1})\right) }}}\\\\
{\frac{{\E^{1/2}\left( {\psi ^2}({\epsilon _1})\right) \int_0^1 {\int_0^t {S\left( s \right)\,ds} \,dW\left( t \right)} }}{{\E\left( \psi '({\epsilon _1})\right) }}}
\end{array}}
\right)
\end{align*}
in probability, where $\Gamma _2$ is defined in \eqref{GD}. This result can be generalized to  non-stationary AR($p$) processes with several real and conjugate complex unit roots.

\end{document}